\documentclass[12 pt]{article}
\usepackage{amsmath,amssymb, braket}
\usepackage[english]{babel}
\usepackage{lipsum}
\usepackage{hyperref}
\date{}
\begin{document}
\title{\textbf{On a relation between $\lambda$-full well-ordered sets and weakly compact cardinals}}
\author{by\\\\Gabriele Gull\`a}
\maketitle
\begin{abstract}
We prove, via transfinite recursion, the existence, inside any linearly ordered set of appropriate regular cardinality $\lambda$,  of a particular kind of well-ordered subsets characterized by the property of \\\textit{$\lambda$-fullness}.
\\Let $\mathfrak{H}$ be a set of regular cardinals: by using our results about well-ordered $\lambda$-full sets we show that if $\inf \mathfrak{H}$ is a weakly compact cardinal, then, for every LOTS $X$, $\mathfrak{H}$-compactness is equivalent to the nonexistence of gaps of types in $\mathfrak{H}$.\\
\end{abstract}
\let\thefootnote\relax
\section*{Introduction}
Our target in this paper is mainly to prove the reverse implication with respect to Corollary 7.3(2) proved in \cite{Lip5}, implication which completes the relation between weakly compact (w.c. so on) cardinals and the notion of \textit{$\mathfrak{H}$-compactness} (see Definition \hyperlink{*}{1.6}) for linearly ordered topological spaces (LOTS).

Let us start by considering the results in \cite{Lip5} characterizing \textit{ultrafilter convergence} and \textit{ultrafilter compactness} inside linearly ordered topological spaces (LOTS) and generalized ordered topological spaces (GOTS).
 Lipparini has produced an extensive literature about the study of covering and converging properties of several kinds of spaces, and the properties of the cartesian product of these spaces (see \cite{Lip1}, \cite{Lip2}, \cite{Lip3} and \cite{Lip4}). These topics are also strongly related to pure logical and set theoretical aspects (see \cite{AH}, \cite{BDM}, \cite{C} and \cite{MS}).\ The mentioned characterizations are shown to be deeply connected with the notions of \textit{gap} and \textit{pseudo-gap} (see Definitions \hyperlink{0}{1.2} and \hyperlink{1}{5.3}).

We recall that given a LOTS $X$, a set $I$ and an ultrafilter $D$ over $I$, a sequence $(x_i)_{i\in I}$ in $X$ is said to $D$-converges to some $x\in X$ if, for every open neighbourhood $U$ of $x$, the set $\left\{i\in I: x_i\in U\right\}$ is an element of $D$.
\\Theorem 3.1 of \cite{Lip5} shows that there is a relation between $D$-convergence and the notion of a gap: if we let $$A:=\left\{x\in X: \left\{i\in I: x<x_i\right\}\in D\right\}$$ and $$B:=\left\{x\in X: \left\{i\in I: x>x_i\right\}\in D\right\},$$ 
then a sequence $(x_i)_{i\in I}$ $D$-converges in $X$ if and only if $(A, B)$ is not a gap.

Let us recall a characterization of the non-existence of gaps in terms of the following notion: a space $X$ is \textit{$D$-compact} if every $I$-indexed sequence in $X$ $D$-converges to some $x\in X$. An ultrafilter $D$ over a set $I$ is \textit{$\lambda$-decomposable}, for a cardinal $\lambda$, if there is a function $f:I\rightarrow \lambda$ such that $f^{-1}(Y)\notin D$ for every $Y\subseteq \lambda$ with $|Y|<\lambda$.
\\Now, let $\mathfrak{K}_D$ be the set of infinite regular cardinals $\kappa$ such that an ultrafilter $D$ is $\kappa$-decomposable.
\\In Theorem 4.1 of \cite{Lip5} it is proved that for every GO-space $X$ (see Definition 5.1) the following are equivalent\\
1) $X$ is $D$-compact\\
2) $X$ is $[\kappa, \kappa]$-compact (see Definition 1.3) for every $\kappa\in \mathfrak{K}_D$\\
3) $X$ has no gaps nor pseudo-gap of type in $\mathfrak{K}_D$\\
4) For every $\kappa\in \mathfrak{K}_D$ and every strictly increasing (resp. decreasing) \\$\kappa$-indexed sequence in $X$, the sequence has a supremum (resp. infimum) to which it converges.

Moreover, for every infinite cardinal $\lambda$ and for an ultrafilter $D$ over $\lambda$ such that there is a family $(Z_{\alpha})_{\alpha\in \lambda}$ of members of $D$ with the property that the intersection of any infinite sub-family is empty, the following holds (Theorem 5.1 of \cite{Lip5}):\\
for 
every GO-space $X$, $X$ is $D$-compact 
if and only if $X$ has no gap nor pseudo-gaps of type $\le\lambda$.

Recall that a regular cardinal $\kappa$ is w.c. if and only if for every LOTS $X$, $[\kappa, \kappa]$-compactness of $X$ is equivalent to the nonexistence of gaps of type $\kappa$ in $X$.
\\We say that a space $X$ is $\mathfrak{H}$-compact if it is $[\kappa, \kappa]$-compact for every $\kappa\in\mathfrak{H}$. 
\\This led to the formulation of the statement that we call $\mathbf{P(\mathfrak{H})}$, where $\mathfrak{H}$ is a class of regular cardinals (\cite{Lip5}):\\
\\\textit{For every LOTS $X$, $X$ is $\mathfrak{H}$-compact if and only if it has no gaps of type in $\mathfrak{H}$}.

In Section 7 of \cite{Lip5} it is shown that if $\mathbf{P(\mathfrak{H})}$ holds, then $\inf \mathfrak{H}$ is a w.c. cardinal (or $\omega$).\\
The purpose of this paper is to prove the converse of Theorem 7.3(2) of \cite{Lip5}: if $\inf\mathfrak{H}$ is a weakly compact cardinal then $\mathbf{P(\mathfrak{H})}$ holds.

Here we are interested in the relation between linearly ordered set with particular well-ordered subsets (namely \textit{$\lambda$-full} ones; for this notion see Definition \hyperlink{2}{2.1}) and linearly ordered topological space with some appropriate covering properties (the mentioned $\mathfrak{H}$-compactness).

In particular, we prove the existence of $\lambda$-full well-ordered sets of various order types (and we prove they have also a stronger version of $\lambda$-fullness).

After the relevant preliminaries (Section 1), in Section 2 we prove that if $\lambda$ is an infinite regular cardinal and $X$ a LO set such that $|X|=\lambda$, then $X$ has a $\lambda^{\leftrightarrows}$-ordered subset or an $\omega^{\leftrightarrows}$-ordered $\lambda$-full subset, and we extend the $\omega^{\leftrightarrows}$-ordered case to every countable $\alpha$ (Theorem \hyperlink{3}{2.2} and Corollary \hyperlink{4}{2.4}).  With the superscript arrows we represent the order ($\alpha^{\rightarrow}$ or just $\alpha$)/reverse order ($\alpha^{\leftarrow}$) of a set.

In Section 3 we prove that if $X$ is a LO set such that $|X|=\lambda>2^{\aleph_{\alpha}}$ for $\lambda$ regular, and if we suppose that 
$$\forall Y\subseteq X\enspace \Bigl(|Y|=\lambda\Rightarrow \exists y\in Y \Bigl(|(-\infty, y)_Y|=|(y, \infty)_Y|=\lambda\Bigr)\Bigr)$$
then $X$ contains an $\omega_{\alpha+1}^{\leftrightarrows}$-ordered $\lambda$-full subset (Theorem \hyperlink{5}{3.2}).

In Section 4 we prove that if $X$ is a LO set of regular cardinality $\lambda\ge \mu$, with $\mu$ w.c., then $X$ admits a $\lambda$-full $\mu^{\leftrightarrows}$-ordered subset or a $\lambda^{\leftrightarrows}$-ordered one (see Theorem \hyperlink{6}{4.6}). This result implies that if $\inf \mathfrak{H}$ is a w.c. cardinal, then $\mathbf{P}(\mathfrak{H})$ holds (Theorem \hyperlink{7}{($\natural$)}).

Finally, in Section 5 we extend the previous results to the GO-spaces case.
\section{Preliminaries}
Here we want to recall the necessary definitions and preliminary results.\\\\
\textbf{Definition 1.1}\\
We say a set $X$ is a \textit{linearly ordered topological space} (LOTS so on) if it is (strictly) linearly ordered (by a relation that we call just $<$) and we consider on it the induced \textit{interval topology}, namely the one generated by the open rays
$$(-\infty, a)=\left\{x\in X| x<a\right\}\quad\mbox{and} \quad (b, \infty)=\left\{x\in X| b<x\right\}$$
for all $a, b\in X$.\\
\\\hypertarget{0}{\textbf{Definition 1.2}}\\
Let $X$ be a LOTS. The ordered pair $(A, B)$ is a \textit{gap} of $X$ if:\\
1) $A$ and $B$ are open subset of $X$;\\
2) $X=A\cup B$;\\
3) $a<b$ for every $a\in A$ and for every $b\in B$;\\
4) $A$ has no maximum and $B$ has no minimum.
\\If $\kappa$ is a cardinal and $(A, B)$ is a gap, we say that $\kappa$ is \textit{a type} of the gap if $cf(A)=\kappa$ or $ci(B)=\kappa$. Thus a gap has at most two types. If a gap has as a type $\kappa$ we say it is a $\kappa$-gap.\\
\newpage\noindent\textbf{Definition 1.3}\\
Let $\kappa\le \lambda$ be cardinals. A topological space is said to be \textit{$[\kappa, \lambda]$-compact} if every open covering of cardinality $\le \lambda$ admits an sub-covering of cardinality $< \kappa$.\\
\\\textbf{Definition 1.4}\\
Let $X$ be a topological space. An accumulation point $x_0\in X$ for an infinite set $A$ is said to be \textit{complete} (c.a.p.) if for every neighbourhood $U$ of $x_0$ we have that
$$|U\cap A|=|A|$$
\\Next theorem is a fundamental bridge between definitions 1.3 and 1.4:\\\\
\textbf{Theorem 1.5} (Alexandroff-Urysohn \cite{AU})\\
Let $\lambda$ be a regular cardinal.\ Then $X$ is $[\lambda, \lambda]$-compact if and only if every subset of cardinality $\lambda$ of $X$ has a c.a.p.\\
\\We will use last characterization often through this work, so sometime we will avoid to recall it.\\
\\Now, let $\mathfrak{H}$ be a class of infinite cardinals.\\
\\\hypertarget{*}{\textbf{Definition 1.6}}\\
A topological space is \textit{$\mathfrak{H}$-compact} if it is $[\kappa, \kappa]$-compact for every $\kappa\in \mathfrak{H}$.\\
\\Now we can give the statement which will be the main motivation for our investigation:\\
\\\hypertarget{ph}{$\mathbf{P}(\mathfrak{H})$} (see \cite{Lip5}): For every LOTS $X$, $X$ is $\mathfrak{H}$-compact $\Leftrightarrow$ $X$ has no gaps of type in $\mathfrak{H}$.
\\\\
\textbf{Remark 1.7}\\
We notice that, if $\lambda$ is regular, a LOTS $X$ $[\lambda, \lambda]$-compact has no gap of type $\lambda$: a strictly increasing cofinal subset of cardinality $\lambda$ in $A$ (or a strictly decreasing coinitial set of the same cardinality in $B$) should have a c.a.p., which contradicts the fact that $(A, B)$ is a gap.

In general $\mathbf{P}(\mathfrak{H})$ is not verified, in fact let us take $X=(\mathbb{R}\times\mathbb{Z}, <_{Lex})$ and take $\mathfrak{H}=\left\{\mathfrak{c}\right\}$. Then $X$ is not $[\mathfrak{c}, \mathfrak{c}]$-compact because it is a space of size $\mathfrak{c}$ with discrete topology.\\ But it has no gap of type $\mathfrak{c}$: for every gap $(A, B)$ the set $(\mathbb{Q}\times \mathbb{Z})\cap A$ is cofinal in $A$ and the set $(\mathbb{Q}\times \mathbb{Z})\cap B$ is coinitial in $B$.
So every gap of $X$ has type at most $\aleph_0$.

More in general, the example works of course also if $\mathfrak{H}=\left\{\lambda\right\}, \omega<\lambda< \mathfrak{c}$ (which is meaningfull only if CH does not hold).
\\But there are cases in which $\mathbf{P}(\mathfrak{H})$ holds: for example if $\mathfrak{H}=Reg\cap[\aleph_0, \lambda]$ (see Corollary 5.1 in \cite{Lip5}, \cite{C}, \cite{GFW}) or if $\mathfrak{H}=\mathfrak{K}_D$ where $mathfrak{K}_D$ is the set of $\kappa$ such that $D$ is $\kappa$-indecomposable, for some ultrafilter $D$ (Corollary 4.1 in \cite{Lip5}).

In this last case, it can be shown (sections 5 and 7 of \cite{Lip2} and section 6 of \cite{Lip5}) that, if GCH holds or if $\mathfrak{H}$ is an interval of cardinals which does not contain the cofinality of its upper extreme, then there is not such an ultrafilter $D$. See Section 7 of \cite{Lip5} for more details.
\section{$\lambda$-Full Sets}
In this section we demonstrate the existence of a particular kind of well-ordered subsets and we explain why it is interesting to have them. 
Let us note that from here to the end of the paper when we use cardinals we assume that they are regular, and $\mathfrak{H}$ will be a class of regular cardinals.\\

First of all we give the following definition, which is a slight modification of that given in section 7 of \cite{Lip5}:\\
\\\hypertarget{2}{\textbf{Definition 2.1}}\\
Let $X$ be a linearly ordered (LO) set of cardinality $\lambda$.
Let $Y\subseteq X$, and let $(y_1, y_2)_X$, with $y_1<y_2\in Y$, be the set
$$\left\{x\in X |\enspace y_1<x<y_2\right\}$$
We say that $Y$ is \textit{$\lambda$-full in X on the right} if
$$\forall y\in Y\qquad \Bigl|\bigcup_{y<y'\in Y}(y, y')_X\Bigr|=\lambda$$
$Y$ is \textit{$\lambda$-full in X on the left} if
$$\forall y\in Y\qquad \Bigl|\bigcup_{Y\ni y'<y}(y', y)_X\Bigr|=\lambda$$
\\
This definition does not depend on the existence of maximum or minimum of $Y$.\
If the context is clear, we will omit the locutions ``on the left''/``on the right''.

We observe that, trivially, if $Y\subseteq X$, with $|X|=\lambda$, is $\lambda$-ordered (resp. $\lambda^{\leftarrow}$-ordered for the reverse order), then it is $\lambda$-full on the right in itself (resp. on the left in itself) into itself assuming $\lambda$ is infinite, so, a fortiori, it is $\lambda$-full in $X$.

We use last definition in the following result. We want to stress the fact that this theorem has been formulated in \cite{Lip5}, but without a proof.
\\\hypertarget{3}{\textbf{Theorem 2.2}}\\
Let $\lambda$ be an infinite regular cardinal.\\ Let $X$ be a LO set such that $|X|=\lambda$.\\
Then $X$ has a $\lambda^{\leftrightarrows}$-ordered subset or an $\omega^{\leftrightarrows}$-ordered $\lambda$-full subset.\\
\\Proof:
\\Case (1): initially let us suppose that 
$$\exists Y\subseteq X\enspace\Bigl(|Y|=\lambda\enspace\wedge\enspace \forall y\in Y \Bigl(|(-\infty, y)_Y|<\lambda\quad\dot{\vee}\quad |(y, \infty)_Y|<\lambda\Bigr)\Bigr)$$
Of course, for each $y\in Y$, it is not possible that both $(-\infty, y)_Y$ and $(y, \infty)_Y$ have size less than $\lambda$.

Using a recursion of lenght $\lambda$ we build a sequence of elements $(y_{\alpha})_{\alpha<\lambda}$ in $Y$ and a sequence of sets $(A_{\alpha})_{\alpha<\lambda}\subseteq Y$ such that, for each $\alpha$, $\enspace y_{\alpha}\in A_{\alpha}$ and the complement $A_{\alpha}^c$ of $A_{\alpha}$ in $Y$ has size $<\lambda$ (so, of course, $|A_{\alpha}|=\lambda$):

-step 0:\\ we choose as $y_0$ any element in $A_0:=Y$;

-successor step:\\ 
now suppose we have already built a sequence $(y_{\beta})_{\beta\le \alpha}$ and a sequence of sets $(A_{\beta})_{\beta\le \alpha}$ as above, and we define $y_{\alpha+1}$ and $A_{\alpha+1}$.\\
We can have two possibilities:\\
\\$a)\quad |(y_{\alpha}, \infty)|=\lambda$
or\\
$b)\quad |(-\infty, y_{\alpha})|=\lambda$\\
\\In the first case we choose any $y_{\alpha+1}\in A_{\alpha+1}:=A_{\alpha}\cap (y_{\alpha}, \infty)$,\\in the second one we choose any $y_{\alpha+1}\in A_{\alpha+1}:=A_{\alpha}\cap (-\infty, y_{\alpha})$.
\\Let us check that, in case (a), $|A_{\alpha+1}^c|<\lambda$: in fact
$$A_{\alpha+1}^c=(A_{\alpha}\cap (y_{\alpha}, \infty))^c=A_{\alpha}^c\cup (-\infty, y_{\alpha})$$
but $|A_{\alpha}^c|<\lambda$ and $|(-\infty, y_{\alpha})|<\lambda$, and so is the same for the union. 
\\Case (b) is,\textit{ mutatis mutandis}, completely analogous.

-limit step $\gamma$: \\
Again, we suppose that $\forall \alpha<\gamma\enspace |A_{\alpha}^c|<\lambda$.\\ 
Then $|\bigcup_{\alpha<\gamma} A_{\alpha}^c|<\lambda$ (because $\lambda$ is regular).\\
Therefore we choose as $y_{\gamma}$ any element in $A_{\gamma}:=Y\setminus \left(\bigcup_{\alpha<\gamma} A_{\alpha}^c\right)$.\\For what we just said, obviously $|A_{\gamma}^c|=|\bigcup_{\alpha<\gamma} A_{\alpha}^c|<\lambda$.
\\We divide the elements of the sequence in the following way by using two subsets $\mathcal{I}$ and $\mathcal{J}$:\\
if $y_{\alpha+1}>y_{\alpha}$ then we put $y_{\alpha}$ in $\mathcal{I}$, while if $y_{\alpha+1}<y_{\alpha}$  we put $y_{\alpha}$ in $\mathcal{J}$.
\\$\mathcal{I}$ and $\mathcal{J}$ define implicitly two disjoint sub-sequences, one strictly increasing and the other strictly decreasing.\\
In fact (increasingness of $\mathcal{I}$): let us suppose $y_{\delta}, y_{\beta}\in \mathcal{I}$ and $\delta>\beta$; if $y_{\beta}\in \mathcal{I}$ this means that $y_{\beta+1}>y_{\beta}$ then $|(y_{\beta}, \infty)|=\lambda$ and then, by construction, $y_{\delta}\in(y_{\beta}, \infty)$, so $y_{\delta}>y_{\beta}$.\\Decreasingness of $\mathcal{J}$ is analogous.
\\So we obtain two sequences such that at least one has length $\lambda$.
\\Clearly we succeed in obtaining a sequence of length $\lambda$ because in both eventualities (a) and (b) we choose always in the set with $\lambda$ elements.

Case (2): now, let us suppose that 
$$\hypertarget{I}{(**)}\quad\forall Y\subseteq X\enspace \Bigl(|Y|=\lambda\Rightarrow \exists y\in Y \Bigl(|(-\infty, y)_Y|=|(y, \infty)_Y|=\lambda\Bigr)\Bigr)$$ 
We build a $\lambda$-full $\omega^{\rightarrow}$-ordered set using the following recursion:

-step 0: we choose an $Y$ of size $\lambda$ (and we call it $Y_0$) and we choose as $y_0$ one of those $y$ (surely existing thank to the hypothesis) such that\\ $|(-\infty, y)_Y|=|(y, \infty)_Y|=\lambda$;

-successor step $n+1$: we consider the set $Y_{n+1}=(y_n, \infty)$ and we choose inside this set $y_{n+1}$ as (one of those) $y$ such that 
$$|(y_{n}, y)|=|(y, \infty)|=\lambda$$
\\So we create a descending chain of sets of size $\lambda$
$$Y=Y_0\supseteq Y_1\supseteq Y_3\supseteq ......\supseteq Y_n\supseteq........$$
which produces, in the way just described, a strictly increasing sequence  $(y_n)_{n\in \omega}$ such that $A=\left\{y_n: n\in \omega\right\}$ is obviously $\omega$-ordered.\\
It is, of course, $\lambda$-full (on the right) by construction.\\
$\mathcal{a}$
\\\\\textbf{Remark 2.3}
\\First we observe that, in the $\lambda^{\leftrightarrows}$-case, because of the choice of $y_{\alpha}$ is free (in the subset of size $\lambda$), in $Y$ we can find many $\lambda^{\leftrightarrows}$-ordered subsets: in fact given one of these, any subset of size $\lambda$ is also $\lambda^{\leftrightarrows}$-ordered.
\\\\In the $\omega^{\leftrightarrows}$-case, the $\lambda$-full sets we built have a stronger property:\\
\\\textbf{Definition $(\bullet)$}\\
We say that a set $Y\subset X$ is \textit{completely $\lambda$-full} if 
$$\forall y, y'\in Y\enspace\mbox{such that}\enspace y<y'\Bigl(|(y, y')_X|=\lambda\Bigr)$$
\\In the proof of Theorem 2.2 we showed the existence of at least one set with the desired properties by building it explicitly.
If, instead of final sets, we had choosen $y_i$'s in initial sets, we would have obtained a strictly decreasing sequence (and so a $\omega^{\leftarrow}$-ordered set) $\lambda$-full on the left.\\
\\Not only, we can find ``a lots'' of other $\omega^{\leftrightarrows}$-ordered $\lambda$-full subsets : starting with an $Y$ of size $\lambda$, we consider the $y$'s which progressively divides $Y$ in sets of size $\lambda$; then it is built a binary tree whose more external branches correspond to the two sequences of orders $\omega$ and $\omega^{\leftarrow}$ respectively in the following way:
$$y_0< y_{12}:=y_1< y_{122}:=y_2< ....$$
and
$$y_0>y_{11}:=y'_1>y_{111}:=y'_2> ....$$
Inside the tree we can then choose branches whose nodes constitute sequences eventually increasing or eventually decreasing, that is sequences which became monotone after a finite number of nodes; from these we can extract sub-sequences of order $\omega$ or $\omega^{\leftarrow}$ erasing initial segments composed by disordered nodes.
\\Or we could be in a case similar to case ``$\lambda$'' with a complete disordered set: in this case we just use the same division with subsets $\mathcal{I}$ and $\mathcal{J}$.\\\\
Now we can extend last part of last result for every countable ordinal:
\\\\\hypertarget{4}{\textbf{Corollary 2.4}}\\
Let $X$ be as in Theorem 2.2, and suppose that $(**)$ holds.\\
Then $X$ contains $\alpha^{\leftrightarrows}$-ordered $\lambda$-full subsets for all countable $\alpha$.
\\\\
Proof\\
We know from the proof of case (2) of Theorem 2.2 that, by using just $(**)$, we can build an $\omega^{\leftrightarrows}$-ordered $\lambda$-full subset inside any subset of size $\lambda$ of $X$.
\\We start by fixing inside $X$ a subset $Y$ of size $\lambda$ and we build a set $\lambda$-full on the right.
\\For the 0 step we build an $\omega$-ordered $\lambda$-full (on the right) set $S_{\omega}$ in $(-\infty, y]_Y$, where $y$ is one of those points which divide $Y$ in two parts of size $\lambda$.
\\For the successor step we do the following: we suppose we can build, by using just $(**)$, a $\lambda$-full $\alpha$-ordered subset in any subset of size $\lambda$ of $X$, and we suppose we built such a $S_{\alpha}$ in $(-\infty, y_{\alpha}]_Y$, where again $y_{\alpha}$ is one of those points which divide $Y$ in two parts of size $\lambda$.\\
Now, $|(y_{\alpha}, \infty)_Y|=\lambda$, so there is an $y'\in (y_{\alpha}, \infty)_Y$ such that $$|(y_{\alpha}, y')_Y|=|(y', \infty)_Y|=\lambda$$ 
So by using $(**)$ we can build a copy $S'_{\alpha}$ of $S_{\alpha}$ inside $(y_{\alpha}, y')_Y$.\\ 
We define $y_{\alpha+1}:=y'$ and 
$$S_{\alpha+1}:=S'_{\alpha}\cup\left\{y_{\alpha+1}\right\}$$
For the limit step $\beta=\sup_{\alpha<\beta} \alpha$ we suppose we have $\alpha$-ordered $\lambda$-full subsets $S_{\alpha}\enspace$ as above, $\forall \alpha<\beta$, so we define
$$S_{\beta}:=\bigcup_{\alpha<\beta}S_{\alpha}$$
We can do this because\\
1) every $S_{\alpha}$ is well-ordered with type $\alpha$ and\\
2) every element of $S_{\alpha}$ is less than every element of $S_{\alpha+1}$, in fact \\$S_{\alpha}\subseteq (-\infty, y_{\alpha}]_Y$ and $S_{\alpha+1}\subseteq (y_{\alpha}, y']_Y$, and of course every element of $(-\infty, y_{\alpha}]_Y$ is smaller than every element of $(y_{\alpha}, y']_Y$.\\
So the union for limit steps guarantees that we obtain a well ordered set.
\\Those subsets are also $\lambda$-full 
by construction (and in fact completely $\lambda$-full).
\\The anti-well-ordered case is analogous by replacing final intervals with initial ones.\\$\mathcal{a}$\\
\\So, from Theorem 2.4 and using the same method of case (1) of the proof of Theorem 2.2, we obtain the following
\\\\
\textbf{Corollary 2.5}\\
Let $\lambda$ be an infinite regular cardinal.\\ Let $X$ be a LO set such that $|X|=\lambda$.\\
Then $X$ has a $\lambda^{\leftrightarrows}$-ordered subset or, for every countable $\alpha$, a $\lambda$-full $\alpha^{\leftrightarrows}$-ordered subset.\\$\mathcal{a}$
\\\\\textbf{Remark 2.6}\\
If $\lambda=\mathfrak{c}$ there are LOTS of cardinality $\lambda$ which do not contain any strictly monotone sequence of cardinality $\kappa$ with $\omega<\kappa\le\mathfrak{c}$ (that is subset $(>\omega)^{\leftrightarrows}$-ordered), and, a fortiori, no one $\lambda$-full.\\
In fact $X=\mathbb{R}$, by using the order density of $\mathbb{Q}$ in $\mathbb{R}$, can not contain any $\omega^{\leftrightarrows}$-ordered subsets.
So, if $\lambda=\mathfrak{c}$, Corollary 2.5 is the best result we can obtain.
\\More in general, a $\kappa$-separable set does not contain a $(>\kappa)^{\leftrightarrows}$-sequence, where we call a set $X$ \textit{$\kappa$-separable} if $X$ contains an order-dense subset of size $\kappa$ (so $\aleph_0$-separability is the usual one).

If a LOTS $X$ has cardinality $>2^{\kappa}$ it cannot be $\kappa$-separable: let us suppose that $X$ contains an order dense subset $B$ of size $\kappa$.\\
Then the function
$$f:X\longrightarrow \mathcal{P}(B)$$
which maps $x$ in $\left\{y\in B: y<x\right\}$, is an injection of $X$ into the power set of $B$, and so $|X|\le 2^{\kappa}$.

Then it follows that, if the LOTS has cardinality $>2^{\kappa}$, it could have a $\kappa^+$-ordered subset.\ We will see this is the case, and we will show that is possible to build well-ordered subsets which are also $\lambda$-full.\\ 
\\Now, separability apart, Remark 2.6 can be generalized by\\
\\\textbf{Proposition 2.7}\\
If $\lambda=2^{\kappa}$ there is at least one LOTS $X$ which does not contain\\ $(>\kappa)^{\leftrightarrows}$-ordered subsets and so no one $\lambda$-full.\\
\\Which follows directly from\\\\ 
\textbf{Lemma 2.8} (Theorem 5.4 of \cite{CN})\\
Let $T$ be the ordered set $\langle\left\{0, 1\right\}^{\kappa}, \le_{lex}\rangle$, with $\kappa\ge \omega$;\\
if $T'\subset T$ is a well (or anti well)-order, it has cardinality at most $\kappa$.\\ 
\\Next theorem is an important connection among $\lambda$-full well-orders, the notion of gap and $[\lambda, \lambda]$-compactness, and it will be crucial for our purpose:
\\\\
\textbf{Theorem 2.9} (Proposition 7.2 of \cite{Lip5})\\
Let $\kappa\le \lambda$ be regular cardinals. The following are equivalent:\\
\\(i) Every LO set of cardinality $\lambda$ has a subset $\lambda^{\leftrightarrows}$-ordered, or a subset\\ $\kappa^{\leftrightarrows}$-ordered $\lambda$-full;\\
(ii) Every LOTS which has no $\kappa$-gap or $\lambda$-gap is $[\lambda, \lambda]$-compact.\\
\section{$\lambda$-Fullness beyond the Countable}
It is quite natural to ask about $\lambda$-full $\omega_{\alpha}^{\leftrightarrows}$-ordered subsets with $\alpha\ge 1$.\\
Here we show they exist for appropriate $\lambda$'s.\\
The existence of well-ordered subsets has been studied quite extensively by Hausdorff and Uryshon, who proved the following:\\
\\\textbf{Theorem 3.1} (Urysohn \cite{U})\\
If $X$ is a LO set of size $>2^{\aleph_{\alpha}}$, then it contains an $\omega_{\alpha+1}^{\leftrightarrows}$-ordered subset.\\
\\Now, it is clear that last theorem does not say anything about the existence of $\lambda$-full well-ordered sets, but we can use this result in our next theorem.
\\\\
\hypertarget{5}{\textbf{Theorem 3.2}}\\
Let $X$ be a LO set such that $|X|=\lambda>2^{\aleph_{\alpha}}$ for $\lambda$ regular, and let us suppose that $(**)$ holds.\\
Then $X$ contains an $\omega_{\alpha+1}^{\leftrightarrows}$-ordered $\lambda$-full subset.\\
\\Proof\\
We start by considering the following equivalence relation on $X$:
$$a\enspace\rho\enspace b\quad\Leftrightarrow\quad \Bigl(a\le b\enspace\wedge\enspace |[a, b]|<\lambda\Bigr)\enspace \vee\enspace\Bigl (b\le a\enspace\wedge\enspace |[b, a]|<\lambda\Bigr)$$
Let us prove that the cardinality of the quotient $X/ \rho$ is $\lambda$: in fact let us suppose there is an equivalence class $H$ with $\lambda$ elements; then for $(**)$ there is a $h\in H$ such that
$$|(-\infty, h)_H|=|(h, \infty)_H|=\lambda$$ 
and so, by $(**)$ again, in $(-\infty, h)_H$ there is a $h'$ such that
$$|(h', h)_H|=\lambda\Rightarrow |[h', h]_X|=\lambda$$
meaning that $h$ and $h'$ are not in relation, contradicting the hypothesis that they belong to the same class.\\
So no one class has $\lambda$ elements and so, by regularity of $\lambda$, the classes must be $\lambda$.

Now we choose a representative $\bar{x_i}$ for every class $K_i$, and we define
$X^*$ as the set of these representatives.
Evidently $X^*\subseteq X$ and it inherits the total order of $X$, so it is a LO set; its cardinality is $\lambda$ and moreover
$$\forall x, y\in X^* (x<y\Rightarrow |(x, y)_X|=\lambda)$$
by construction.

For Theorem 3.1, $X^*$ admits an $\omega_{\alpha+1}^{\leftrightarrows}$-ordered subset $Z$; now, every element of $Z$ is an element of $X^*$, so for every couple of elements $z_i<z_j\in Z$ we have that
$$|(z_i, z_j)_X|=\lambda$$
and so $Z$ is completely $\lambda$-full and then $\lambda$-full.\\$\mathcal{a}$\\
\\\textbf{Remark 3.3}\\
For what we just proved, it follows that, under $(**)$, every subset of size $\lambda$ of $X$ contains an $\omega_{\alpha+1}^{\leftrightarrows}$-ordered $\lambda$-full subset, so, by using the same recursion of the proof of Corollary 2.4, $X$ contains a $\lambda$-full $\gamma^{\leftrightarrows}$-ordered subset for every $\gamma$ of cardinality $\aleph_{\alpha+1}$: we start with an $\omega_{\alpha+1}^{\leftrightarrows}$-ordered $\lambda$-full subset and, by using hypothesis $(**)$ as in Corollary 2.4, we build bigger (from the point of view of the order) $\lambda$-full subsets.
\\\\
The previous theorem together with the argument for case (1) of the proof of Theorem 2.2 yields the following\\\\
\textbf{Corollary 3.4}\\
Let $X$ be a LO set s.t. $|X|=\lambda>2^{\aleph_{\alpha}}$ for $\lambda$ regular.\\
Then $X$ admits an $\omega_{\alpha+1}^{\leftrightarrows}$-ordered $\lambda$-full subset  or a $\lambda^{\leftrightarrows}$-ordered one.\\
\section{Weakly Compact Cardinals and $\lambda$-Full Well-Ordered Sets}
In this section we state a connection among w.c. cardinals, $\lambda$-full well-ordered sets and $\mathfrak{H}$-compactness of LOTS.\\
We remember that a cardinal $\kappa$ is w.c. iff every linear order of cardinality $\kappa$ has a $\kappa^{\leftrightarrows}$-ordered subset.\\
We notice that, thanks to Theorem 2.2, every LO set of cardinality $\aleph_0$ has an $\omega^{\leftrightarrows}$-ordered subset.
\\Next theorem is the first link between w.c. cardinals and property $\mathbf{P}(\mathfrak{H})$, and is a sort of root for our last result:\\
\\\textbf{Theorem 4.4} (Corollary 7.3 of \cite{Lip5})\\
Let $\mathfrak{H}$ be a class of regular cardinals. The following hold:\\
1) if $\mathbf{P}(\mathfrak{H})$ holds, then $\inf \mathfrak{H}$ is $\omega$ or a weakly compact cardinal.\\
2) if $\inf \mathfrak{H}=\omega$ then $\mathbf{P}(\mathfrak{H})$ holds.\\
\\
With these elements we want to prove the following statement:
\newpage\noindent\hypertarget{7}{\textbf{Theorem ($\natural$)}}\\ 
\textit{If $\inf \mathfrak{H}$ is a w.c. cardinal, then $\mathbf{P}(\mathfrak{H})$ holds}.
\\\\ 
A proof of ($\natural$), together with Theorem 4.4, would characterize completely $\mathbf{P}(\mathfrak{H})$:\\
\\\textbf{Theorem $(\flat)$}\\ 
\textit{$\mathbf{P}(\mathfrak{H})$ holds if and only if $\inf \mathfrak{H}$ is $\omega$ or a w.c. cardinal}.\\
\\In order to prove Theorem ($\natural$) we first prove the following\\
\\\hypertarget{6}{\textbf{Theorem 4.6}}\\
If $X$ is a LO set of regular cardinality $\lambda\ge \mu$, with $\mu$ w.c., then $X$ admits a $\lambda$-full $\mu^{\leftrightarrows}$-ordered subset or a $\lambda^{\leftrightarrows}$-ordered one.\\
\\Proof\\
If $\mu$ is w.c., from Lemma 4.3 we know that every LO of size bigger or equal than $\mu$ w.c. has a subset $\mu^{\leftrightarrows}$-ordered.

If we are in case (1) of the proof of Theorem 2.2 and we obtain a $\lambda^{\leftrightarrows}$-ordered set (which is a $\mu^{\leftrightarrows}$-ordered one if $\lambda=\mu$).

Otherwise, we are in case (2) of the proof of Theorem 2.2, then with the same construction used in the proof of Theorem 3.2 we obtain a subset still $\mu^{\leftrightarrows}$-ordered and completely $\lambda$-full, so $\lambda$-full.\\
$\mathcal{a}$
\\\\Proof of Theorem ($\natural$)\\
We have to prove that, if $\inf \mathfrak{H}$ is a w.c. cardinal, then for every LOTS $X$, $X$ is $\mathfrak{H}$-compact iff $X$ has no gaps of type in $\mathfrak{H}$.\\
Theorem 4.6 says that condition (i) of Theorem 2.9 is satisfied when $\kappa=\mu$ w.c., exactly as Theorem 2.2 implied it for $\kappa=\omega$.\\ 
So, by applying (i) as in proof of Theorem 4.4(2) and thanks to Remark 1.7, it follows the thesis.\\
$\mathcal{a}$
\section{The GO-spaces Case}
Here we want briefly extend the results of last section to GO-spaces.\\
\\\textbf{Definition 5.1}\\
A \textit{generalized ordered space} (GO-space or GOTS) is a linearly ordered set with a $T_2$ topology having a base of order-convex sets.\\
\\Examples of GO-spaces are the space of all countable ordinals $[0, \omega_1)$ with the natural order topology, $\mathbb{R}$ with the topology generated by the sets $[a, b)$ (the \textit{Sorgenfrey line}) and $\mathbb{R}$ with the topology generated by the sets \\$\left\{U\cup K: \mbox{$U$ open in $\mathbb{R}$ and $K\subset \mathbb{R\setminus \mathbb{Q}}$}\right\}$ (the \textit{Michael line}).

\v{C}ech showed that GO-spaces are precisely the subspaces of LOTS, and moreover that there is a canonical construction that produces, for any GO-space $X$, a LOTS $X^*$ that contains X as a closed subspace.\\\\
\textbf{Remark 5.2} (see \cite{BL}, \cite{L}, \cite{P3})\\
It is not unusual that a GO-space $X$ has a property if and only if the LOTS $X^*$ has the same property - for example metrizability and Lindel\"of property - but not necessary: the Sorgenfrey line $S$ is a \textit{perfect space} (every closed subset of $S$ is the countable intersection of open sets), but $S^*$ is not.

Anyway, the space $(\mathbb{R}\times \left\{0, 1\right\}, <_{Lex})$ is a perfect LOTS which contains $S$ as a subspace.\\
It happens also that statements involving LOTS are equivalent to statements involving GO-spaces. See \cite{BL} for examples about this.

A particular example, for what concerns us, is a characterization of w.c. cardinals via the notion of LOTS: it is equivalent to the one with GO-spaces with no gaps nor \textit{pseudo-gaps} (see next definition) of type $\kappa$, in place of LOTS: a regular cardinal $\kappa>\omega$ is weakly compact if and only if for every GOTS (equivalently LOTS) $[\kappa, \kappa]$-compactness is equivalent to having no gaps nor pseudo-gaps of type $\kappa$.\\
\\\hypertarget{1}{\textbf{Definition 5.3}}\\
An ordered pair $(A, B)$, where $A$ and $B$ are open subsets of a GO-space $X$, both non-empty, with $A\cup B=X$ and $a<b$ for every $a\in A, b\in B$, is a \textit{pseudo-gap} for $X$ if either $A$ has a maximum and $B$ has no minimum, or $A$ has no maximum and $B$ has a minimum.\\
If $B$ has a minimum the \textit{type} of $(A, B)$ is the cofinality of $A$, while if $A$ has a maximum the type is the coinitiality of $B$.\\
A pseudo-gap can occur just in GO-spaces that are not LOTS because $A$ and $B$ must be open sets.\\
\\Now, Theorem 2.2 holds also for GO-spaces because about $X$ we used just its linear order and its cardinality.\\
In Theorem 2.9, by considering a GO-space in place of a LOTS, and by adding  the non-existence of pseudo-gaps, we have that:\\\\$(i) \Leftrightarrow (\sharp$: For every GO-space $X$, if $X$ has no gaps nor pseudo-gap of type $\kappa$ or $\lambda$, then it is $[\lambda, \lambda]$-compact).\\\\
The proof follows from $(i)\Rightarrow\sharp\Rightarrow (ii)\Rightarrow (i)$, where the first and the third implications are proved in 2.9, and the second one is trivial (by definitions of GOTS and LOTS).\\\\
So we can replace LOTS with GO-space (by adding nonexistence of pseudo-gaps of the same type too) also in 4.4(2).\\\\
We introduce the \textit{generalized} $\mathbf{P}(\mathfrak{H})$\\
\\$\mathbf{GP}(\mathfrak{H})$: For every GO-space $X$, $X$ is $\mathfrak{H}$-compact $\Leftrightarrow$ $X$ has no gaps nor pseudo-gap of type in $\mathfrak{H}$\\
\\and we conclude, exactly as for Theorem ($\natural$), that\\\\ 
\centerline{$\inf\mathfrak{H}$ is a w.c. cardinal (or $\omega$) $\Leftrightarrow$ $\mathbf{GP}(\mathfrak{H})$ holds.}\\\\
So we can summarize the results obtained in the following theorem:\\
\\\textbf{Theorem 5.4}\\
Let $\mathfrak{H}$ be a class of infinite regular cardinals such that $\omega<\kappa=\inf \mathfrak{H}$ and let $X$ be a GO-space (equivalently LOTS).\ Then the following are equivalent:\\
a) $\kappa$ is a weakly compat cardinal\\
b) $X$ is $\mathfrak{H}$-compact $\Leftrightarrow$ $X$ has no gaps nor pseudo-gap of type in $\mathfrak{H}$\\$\mathcal{a}$\\\\
\\The characterization recalled in Remark 5.2 can be seen as a consequence of 5.4 by choosing $\mathfrak{H}=\left\{\kappa\right\}$.

\end{document}